\documentclass{amsart}

\usepackage{amsmath,amssymb,amsthm,latexsym}
\usepackage{hyperref}

\theoremstyle{plain}
\newtheorem{theorem}{Theorem}[section]
\newtheorem{lemma}[theorem]{Lemma}
\newtheorem{proposition}[theorem]{Proposition}
\newtheorem{fact}[theorem]{Fact}
\newtheorem{corollary}[theorem]{Corollary}

\theoremstyle{definition}
\newtheorem{definition}[theorem]{Definition}
\newtheorem{remark}[theorem]{Remark}
\newtheorem{conjecture}[theorem]{Conjecture}

\newcommand{\dc}{\mathrm{dc}}
\newcommand{\semi}{\rtimes}
\newcommand{\len}{\unlhd}
\newcommand{\Q}{\mathbb Q}
\newcommand{\F}{\mathbb F}
\newcommand{\car}{\mbox{char}}
\newcommand{\N}{\mathbb N}

\begin{document}
	
	\title{Sharply $2$-transitive groups of finite Morley rank}          
	\author{Tuna Alt\i nel, Ay\c se Berkman and Frank O. Wagner}
	\address{Université Lyon 1, Centrale Lyon, INSA Lyon, Université Jean Monnet, CNRS, ICJ UMR5208, 69622 Villeurbanne, France}
	\email{altinel@math.univ-lyon1.fr}
	
	\address{Mimar Sinan G\"uzel Sanatlar \"Universitesi;  MSGS\"U Matematik B\"ol\"um\"u, \c Silah\c s\"or Cad. 71, Bomonti, \c Si\c sli, Istanbul 34380, Turkey}
	\email{ayse.berkman@msgsu.edu.tr}
	
	\address{Université Lyon 1, Centrale Lyon, INSA Lyon, Université Jean Monnet, CNRS, ICJ UMR5208, 69622 Villeurbanne, France}
	\email{wagner@math.univ-lyon1.fr}
	\thanks{The first and third authors would like to acknowledge the support of ANR-13-BS01-0006 ValCoMo and AAPG2019 (ANR-DFG) GeoMod. The third author is partially supported by the Science Committee of the Ministry of Science and Higher Education of the Republic of Kazakhstan (Grant No. AP19677451).}
	\date{2 June 2026}
	
	
\begin{abstract}
A sharply $2$-transitive permutation group of finite Morley rank and characteristic $2$ splits; a split sharply $2$-transitive permutation group of finite Morley rank and characteristic different from $2$ is the group of affine transformations of an algebraically closed field. In particular, a sharply $2$-transitive permutation group of finite Morley rank of characteristic $3$ is the group of affine transformations of an algebraically closed field of characteristic $3$.

Without any assumption on Morley rank, a sharply $2$-transitive permutation group of characteristic $0$ splits if its point stabilizers are virtually abelian.
\end{abstract}

\maketitle
\noindent{\em Dedicated to Katrin Tent on the occasion of her 60th birthday}

\section{Introduction}
A {\em Frobenius group} is a group $G$ together with a malnormal subgroup $B$, i.e.\ $B \cap B^g=\{1\}$ for all $g\in G\setminus B$; the group $B$ is called the {\em Frobenius complement}. A Frobenius group $G$ with Frobenius complement $B$ {\em splits} if there is a normal subgroup $N$, called the {\em Frobenius kernel}, such that $G=N\semi B$. All finite Frobenius groups split \cite{Fr01}; moreover the Frobenius kernel is nilpotent \cite{Th60}. This can be seen as a precursor to the Feit-Thompson theorem on the solubility of groups of odd order, which in turn is the starting point for the classification of the finite simple groups.  As Tao has remarked, all known proofs of Frobenius' Theorem use group characters; removing character theory might lead to a new route to the classification of finite simple groups \cite{tao}.

On the model-theoretic side, the classification of the finite simple groups has inspired Borovik's programme for the resolution of
the {\em Algebraicity Conjecture} by Cherlin and Zilber, which asserts that a simple group of finite Morley rank should be an algebraic group over an algebraically closed field. This has been successful in case there is an infinite elementary abelian $2$-subgroup \cite{ABC}. However, since neither representation nor character theory are available in this context, there is no analogue of the Feit-Thompson Theorem, and there may well exist simple groups of finite Morley rank without involutions (the so-called {\em degenerate} case), which are impervious to Borovik's approach.

The interest in the study of Frobenius groups of finite Morley rank is thus two-fold: On one hand, following Tao, one hopes to reproduce Frobenius' Theorem in this context, and to climb up the ladder to eventually eliminate degenerate simple groups of finite Morley rank (i.e.\ prove a finite Morley rank version of the Feit-Thompson Theorem); on the other hand one seeks to develop methods in the finite Morley rank context which might help to find a character-free proof of Frobenius' Theorem for finite groups.

Frobenius groups arise naturally as transitive non-regular permutation groups such that only the identity fixes two distinct points, the Frobenius complement being the stabiliser of a single point. Borovik and Nesin have conjectured that a Frobenius group of finite Morley rank has a nilpotent Frobenius kernel (and in particular splits), see also \cite{CT2}.

A particular case of Frobenius groups are sharply $2$-transitive permutation groups. It is easy to see that they contain an involution: any permutation $g\in G$ exchanging two points $x$ and $y$ must have order $2$; moreover $g$ is the only permutation exchanging $x$ and $y$ by sharp $2$-transitivity, and if $g'$ is an involution exchanging $x'$ and $y'$, then $g'=g^h$ for the unique $h\in G$ with $h(x')=x$ and $h(y')=y$. Thus all involutions are conjugate. Moreover, if the group is split, then the Frobenius kernel is abelian~\cite{Neu}.

It was a long-standing open question whether all infinite sharply $2$-transitive groups split or not. In 2017 Rips, Segev and Tent constructed the first family of examples of non-split sharply $2$-transitive groups \cite{RST} (just called {\em examples} for the rest of this section). After that important result, the interest in this area grew, and many more examples were constructed. In particular, Tent and Ziegler gave a simpler and shorter construction in \cite{TZ}. All of these examples are of permutation characteristic $2$, i.e.\ involutions have no fixed points.

When involutions have fixed points in a sharply $2$-transitive group, the products of two distinct involutions have the same odd prime order $p$, or they are all torsion-free. We say that the (permutation) characteristic is $p$ in the former case, and $0$ in the latter. It is now known that there are examples in characteristic $0$, see \cite{RT,ST}, and also in characteristic $p$, for any  sufficiently large prime number $p$, see \cite{AAT}. Linear examples were constructed in \cite{GG}; examples of simple groups  first appeared in \cite{AT}, and examples of finitely generated simple groups were later given in \cite{AG}. 

On the other hand, one may search for conditions when a sharply $2$-transitive group is necessarily split. An old result by Kerby and Wefelscheid \cite{KW} shows that sharply $2$-transitive groups of characteristic $3$ split (see \cite{tuek} for a short proof). Under some characteristic restrictions, linear sharply $2$-transitive groups split \cite{GG2}; we shall show that sharply $2$-transitive groups of characteristic $0$ whose point stabilizers are virtually abelian split. This answers positively question 12.48(b) of the Kourovka notebook \cite{Kh}.

Currently, no non-split example is known in the finite Morley rank context. We shall call a sharply $2$-transitive groups of finite Morley rank {\em standard} if it is isomorphic to $\operatorname{AGL_1}(K)$ for some (algebraically closed) field $K$. Borovik and Nesin conjectured:
\begin{conjecture}\label{Conjecture} \cite{BN} An infinite sharply $2$-transitive permutation group of finite Morley rank is standard.
More precisely:\begin{enumerate}
\item[(i)] A sharply $2$-transitive permutation group of finite Morley rank splits.
\item[(ii)] A sharply $2$-transitive split permutation group of finite Morley rank is standard.
\end{enumerate}\end{conjecture}
We shall show (i) in characteristic $2$, and (ii) in characteristic different from $2$
(in characteristic $0$ this was already done by Cherlin, Grund\-h\"ofer, Nesin and V\"olklein \cite{CGNV}). In particular, Conjecture \ref{Conjecture} holds in characteristic $3$.
Note that an infinite sharply $2$-transitive group of finite Morley rank is connected, and so is its Frobenius complement~\cite{BN}, i.e.\ they have no definable subgroups of finite index.

However, our results do not advance the search for a finite Morley rank version of the Feit-Thompson Theorem, nor address the question of degenerate simple groups of finite Morley rank, since they are based on the study of involutions whose existence is guaranteed by sharp $2$-transitivity. They might, however, contribute to the study of the {\em odd characteristic} case of the Algebraicity Conjecture (where the group contains a Pr\"ufer $2$-subgroup and has finite $2$-rank), which is still open.

For background material on groups of finite Morley rank the reader should consult \cite{BN}, in particular Chapters 11.3 and 11.4 on Frobenius and sharply $2$-transitive permutation groups.
The results of Section \ref{sec2} were obtained by the third author and appeared online as \cite{W}.

Finally, we should like to thank the anonymous referee for his/her careful reading of the paper.

\section{Near-domains and near-fields.}\label{sec2}
V. D. Mazurov asked in the Kourovka Notebook \cite[question 12.48]{Kh}:\\
Let G be a sharply $2$-transitive permutation group.\begin{enumerate}
\item[(i)] Does G possess a regular normal subgroup if a point stabilizer is locally finite?
\item[(ii)] Does G possess a regular normal subgroup if a point stabilizer has an abelian
subgroup of finite index?\end{enumerate}
We shall answer question (ii) affirmatively in characteristic $0$. We shall work in the equivalent setting of near-domains.

\begin{definition}
$(K,0,1,+,\cdot)$ is a {\em near-domain} if for all $a,b,c\in K$ we have:
\begin{enumerate}
\item $(K,0,+)$ is a {\em loop}, i.e.\ $a+x=b$ and $y+a=b$ have unique solutions, with $a+0=0+a=a$;
\item there is a unique additive two-sided inverse, i.e.\ $a+b=0$ iff $b+a=0$;
\item $(K\setminus\{0\},1,\cdot)$ is a group, and $0\cdot a=a\cdot 0=0$;
\item left distributivity holds: $a\cdot(b+c)=a\cdot b+a\cdot c$;
\item for all $a,b\in K$ there is $d_{a,b}\in K$ such that $a+(b+x)=(a+b)+d_{a,b}\cdot x$ for all $x$.\end{enumerate}
A near-domain is a {\em near-field} if in addition it is associative, i.e.\ $d_{a,b}=1$ for all $a,b$.\end{definition}
Hence a near-field is a skew field iff right distributivity holds. Note that in a near-domain $an=a+\cdots+a$ ($n$ times) by left distributivity, where $n=1+\cdots+1$, but we may have $na\not=an$. 

\begin{fact}[Tits, Karzel \cite{Ka}]\label{f:nearfield}
A sharply $2$-transitive permutation group $G$ is iso\-morphic to the group of affine transformations of some near-domain $K$, i.e.\ to the set of permutations $\{x\mapsto a+bx:a,b\in K,\,b\not=0\}$; in characteristic $\not=2$ the centraliser of any involution is isomorphic to the multiplicative group $K^\times$. The group $G$ is split iff $K$ is a near-field.\end{fact}

The following fact shall be used frequently without further reference.

\begin{fact}[{\cite[(8.1), (8.3), (8.10), (8.16), (8.18)]{Ker74}}]\label{fact} For all $a,b,c\in K$ we have:\begin{enumerate}
\item $d_{a,an}=1$ for all $n\in\N$ (power-associativity of the additive loop).
\item\label{commute} $d_{a,b}(b+a)=a+b$ (so $d_{a,b}d_{b,a}=1$, and $a$ and $b$ commute iff $d_{a,b}=1$).
\item\label{conjugate} $cd_{a,b}c^{-1}=d_{ca,cb}$.
\item\label{addc} $d_{a,b}=d_{a,c}d_{c+a,-c+b}d_{-c,b}$.
\item $|K^\times:C_{K^\times}(d_{a,b})|=\infty$ if $d_{a,b}\not=1$.
\item If $a,b\in E$ then $(a+b)\,2\in E$.
\item In characteristic $>2$ there is a unique maximal sub-near-field, namely $E$.
\end{enumerate}
\end{fact}

Let $E$ be the set $\{d\in K:1+d=d+1\}$. Then $E=\{a\in K:d_{a,1}=1\}$ by Fact \ref{fact}(\ref{commute}). By power-associativity, $1$ generates a subfield of $K$ contained in $E$, which is either $\Q$ or $\F_p$. Thus $K$ has a characteristic, which is easily seen to be equal to the permutation characteristic of $G$.

Let now $A$ be any subgroup of finite index in $K^\times$ which avoids all non-trivial coefficients $d_{a,b}$ for $a,b\in K$.
Kerby \cite[Theorem 8.26]{Ker74} has shown that $K$ must be a near-field in the following cases:\begin{enumerate}
\item $\car(K)=0$ and $|K^\times:A|=2$,
\item $\car(K)=2$, $|K^\times:A|=2$ and $|E|>2$,
\item $\car(K)=p>2$ and $|K^\times:A|<|E|$.
\end{enumerate}
We shall adapt the proof of (3) to characteristic $0$.

\begin{lemma}\label{lemma} Suppose $a,b\in K$ with $d_{a,b}=1$. Then $d_{a,bn}=1$ for all $n\in\N$.\end{lemma}
\begin{proof} E use induction on $n$. This is clear for $n=0$ and $n=1$. So suppose it holds for $n$. Then $d_{a,b}=1$, $d_{a,bn}=1$ and $d_{b,bn}=1$, whence
\[\begin{aligned}b(n+1) +a&=(bn+b)+a= bn+(b+a)= bn+(a+b)\\
&=(bn+a)+b=(a+bn)+b=a+(bn+b)=a+b(n+1).\qedhere\end{aligned}\]
\end{proof}

\begin{proposition}\label{propn} If $A\le K^\times$ is a subgroup of finite index avoiding all nontrivial $d_{a,b}$ and $\car(K)=0$, then $K$ is a near-field.\end{proposition}
\begin{proof} Recall that $\Q\subseteq E$. If $K=E$, then $d_{a,b}=bd_{b^{-1}a,1}b^{-1}=1$ for all $a,b\in K^\times$ and $K$ is a near-field. So assume that $E\subsetneq K^\times$, and take $a\in K\setminus E\,2^{-1}$. Let $n=|K^\times:A|$. Then there are distinct $i>j$ in $\{0,1,2,\ldots,n\}/n!$ with $d_{a,i}A=d_{a,j}A$; since $d_{-j,i}=1$ we obtain by Fact \ref{fact}(\ref{addc}) that
\[d_{a,i}=d_{a,j}d_{j+a,-j+i}d_{-j,i}=d_{a,j}d_{j+a,-j+i}.\]
Hence $d_{j+a,-j+i}\in A$, and $d_{j+a,-j+i}=1$ by assumption. By Fact \ref{fact}(\ref{conjugate}) we have
\[d_{(i-j)^{-1}(j+a),1}=(i-j)^{-1}d_{j+a,-j+i}(i-j)=1,\]
so $(i-j)^{-1}(j+a)\in E$. Since $-(i-j)^{-1}j\in\Q\subseteq E$, we have 
\[[-(i-j)^{-1}j+(i-j)^{-1}(j+a)]\,2=(i-j)^{-1}a\,2\in E,\]
and $d_{a2,i-j}=1$. But $0<(i-j)\,n!\le n$ is integer, and there is an integer $k>0$ with $i-j=\frac1k$. By Lemma \ref{lemma} we obtain $d_{a2,1}=1$ and $a\,2\in E$, a contradiction.
\end{proof}

\begin{corollary} Let G be a sharply $2$-transitive permutation group of characteristic $0$ whose point stabilizer is virtually abelian. Then $G$ is split.\end{corollary}
\begin{proof} If $K$ is the associated near-domain, $K^\times$ has an abelian subgroup $A$ of finite index. Now any non-trivial $d_{a,b}$ has a centralizer of infinite index in $K^\times$, so $d_{a,b}\notin A$. We finish by Proposition \ref{propn}.\end{proof}

\section{Characteristic $2$}
From now on, we assume that the structures considered have finite Morley rank.

In \cite{BN} Borovik and Nesin study mostly Frobenius groups whose Frobenius complement contains an involution. The following result elucidates what happens otherwise. It relies on the even case of the Algebraicity Conjecture.

Recall that a definable subgroup $H\le G$ is {\em generous} if $\bigcup_{g\in G}H^g$ is generic in $G$. If $G$ is connected, any two generic subsets have generic intersection. Hence any two generous subgroups can be conjugated to intersect non-trivially. Examples of generous subgroups include a malnormal subgroup $B$ (as $RM(\bigcup_{g\in G}B^g)=RM(G/B)+ RM(B)=RM(G)$) and 
the centralizer of a {\em decent torus}, i.e.\ a definable divisible abelian subgroup with dense torsion \cite{Ch05}. 
\begin{theorem}\label{Frobenius} Let $G$ be a connected Frobenius group of finite Morley rank with Frobenius complement $B$. If $B$ does not contain an involution, then $G$ has a normal definable connected subgroup $N$ containing all involutions, such that $N\cap B=\{1\}$.\end{theorem}
\begin{proof} If $G$ has no involutions, we may simply take $N=\{1\}$. So we may assume that $G$ has involutions. Note that $B$ is definable by \cite[Proposition 11.19]{BN}. By \cite{BBC} a $2$-Sylow subgroup $S$ of $G$ is infinite. If $S$ has a non-trivial abelian divisible subgroup $T$, then its {\em definable hull} $\dc(T)$, the smallest definable group containing $T$, is a decent torus, its centralizer $C_G(T)=C_G(\dc(T))$ is generous, and after conjugation intersects $B$ non-trivially. But then $T\le B$ by malnormality, a contradiction. It follows that $S$ has bounded exponent.

Let $N$ be a definable connected normal subgroup of $G$ with $B\cap N=\{1\}$; we may choose it of maximal Morley rank possible. Then either $G=N\semi B$ splits or $BN/N$ is malnormal in $G/N$ by \cite[Lemma 11.37]{BN}. Note that $BN/N$ does not contain an involution, so we are done in the first case; in the second case either $N$ contains all involutions and we are done, or we may divide out by $N$ and assume that $G$ has no definable infinite normal subgroup intersecting $B$ trivially. If now $N\len G$ is finite with $B\cap N=\{1\}$, then $N\le Z(G)$ by connectivity of $G$, so $N\le B$ by malnormality, and $N=\{1\}$.

It follows that $G$ has no non-trivial definable normal subgroup intersecting $B$ trivially. In particular it has no non-trivial definable abelian normal subgroup by malnormality of $B$. Hence $G$ has a definable infinite simple normal subgroup $N$ (which may be $G$ itself). Then $N\cap B$ is malnormal, whence generous in $N$, as is $(N\cap B)^g$ for any $g\in G$. It follows that there is $n\in N$ with $(N\cap B)\cap(N\cap B)^{gn}\not=\{1\}$. Then $gn\in B$ by malnormality, and $G=BN$. Moreover, $N$ contains an involution, and is algebraic by \cite{ABC}. But an algebraic Frobenius group splits \cite[Lemma 11.39]{BN}, contradicting simplicity of $N$.\end{proof}
\begin{remark} A simple $B$-invariant section of $N$ does not contain an involution. If it did, it would be algebraic by \cite{ABC}, but a simple algebraic group cannot be a Frobenius complement \cite[Corollary 11.40]{BN}.\end{remark} 

\begin{corollary} An infinite sharply $2$-transitive permutation group of finite Morley rank and characteristic $2$ splits.\end{corollary}
\begin{proof} Let $G$ be the group, $B$ its Frobenius complement, and $N\len G$ the definable normal subgroup given by Theorem \ref{Frobenius}. If $i\in N$ is an involution, then $RM(N)\ge RM(i^B)=RM(B)$. Thus $2RM(B)\le RM(BN)\le RM(G)=2RM(B)$. It follows that $G=N\semi B$ splits.\end{proof}

The splitting result in characteristic $2$ allows to extend the classification of Delahan and Nesin of sharply $2$-transitive groups of finite Morley rank of characteristic different from $2$ with nilpotent point stabilizers (see \cite[Corollary 11.74]{BN}) to the missing characteristic as well, assuming only solubility of the point stabilizer:

\begin{corollary} An infinite sharply $2$-transitive permutation group of finite Morley rank and characteristic $2$ with soluble point stabilizers is isomorphic to $K_+\semi K^*$ where $K$ is an algebraically closed field of characteristic $2$.\end{corollary}
\begin{proof} Let $G$ be a sharply $2$-transitive group as in the statement. By the previous corollary and an application of \cite[Lemma 11.46]{BN}, $G=A\semi H$ where $A$ is an infinite elementary abelian $2$-subgroup while $H$ is a point stabilizer. Since $H$ is assumed to be soluble, $G$ is solvable. One concludes using \cite[Proposition 11.62]{BN}.\end{proof}

\section{Centrality of the Sylow $2$-Subgroup}
In this section we shall give a criterion when the Sylow $2$-subgroup must be central. It will be the main ingredient in Section \ref{five} to show that a split sharply $2$-transitive permutation group in characteristic different from $2$ is standard.
\begin{theorem}\label{t:8} Let $G$ be a connected group of finite Morley rank whose connected definable abelian subgroups are decent tori. Then its $2$-Sylow subgroup is connected and central.\end{theorem}
%, and maximal tori are self-normalizing.\end{theorem}
Note that the hypotheses actually imply that connected definable abelian subgroups are good tori.
\begin{proof} If $A$ is a decent torus, then $C_G(A)$ is generous in $G$ and $A$ does not have an infinite family of definable subgroups \cite{Ch05}. Moreover, $N_G(A)^0=C_G(A)^0=C_G(A)$ by \cite{AB08}.

We claim that a connected soluble subgroup $S$ is abelian: If $A$ is the last non-trivial derived group, then $A$ is central in $S$. So if $N$ is the second last derived subgroup, then $N$ is connected nilpotent, and it is the central product of a divisible abelian group with a connected nilpotent group of bounded exponent. If the latter were infinite, it would contain an infinite elementary abelian subgroup, a contradiction. So $N$ is abelian, and $N=A=S$.

We shall call a maximal definable connected abelian subgroup of $G$ a {\em full torus}; if $A$ is a full torus, then $A=C_G(A)$, and $A$ is almost self-normalizing and generous in~$G$. Recall that by \cite{Ch05} all full tori in $G$ are conjugate.

A connected definable subgroup $H$ of $G$ is called {\em full} if it contains a full torus. We claim that for full $H$ every element $g\in H$ is contained in a full torus $T\le H$. So consider a non-trivial element $g\in H$. Then $C_H(g)$ is infinite \cite{BBC} and contains a decent torus $T$, which we take maximal possible. So $C_H(T)$ is connected, and $g\in C_H(T)$. If $g\notin T$ then $gT$ has an infinite centralizer $C/T$ in $C_H(T)/T$, which contains an infinite connected abelian subgroup $A/T$. But then $A$ is again abelian connected, whence a torus, and
\[g\in C_{C_H(T)}(A/T)=C_{C_H(T)}(A/T)^0\le N_{C_H(T)}(A)^0=C_H(A).\]
Thus $T<A\le C_H(g)$, a contradiction. It follows that every element is contained in a torus $T$. But now $T$ can be extended to a maximal connected abelian subgroup $T'$, which must be conjugate to any full torus contained in $H$ and hence a full torus itself. In particular for $H=G$ we see that any $g\in G$ is contained in some full torus $T\le C_G(g)^0$; as full tori are conjugate, $g$ is contained in all full tori of $C_G(g)^0$.

Let $S$ be the connected component of a Sylow $2$-subgroup (in short, a {\em $2$-Sylow$^0$}). Suppose for a contradiction that $S$ is not central. As $C_G(S)^0=N_G(S)^0$, there is $a\in G$ with $S^a\not=S$. Choose two full tori $A$ and $A'$ containing distinct $2$-Sylows$^0$ $S\not=S'$ such that $A\cap A'$ is maximal possible (full tori are conjugate and hence uniformly definable, as is the intersection of two of them, and there must be a maximal one by stability). Replacing $G$ by $C_G(A\cap A')^0/(A\cap A')$, we may assume that this intersection is trivial. In particular any two full tori with non-trivial intersection have the same $2$-Sylow$^0$. So for any non-trivial $g\in G$, since any two full tori in $C_G(g)^0$ are conjugate and must contain $g$, there is a unique $2$-Sylow$^0$ in $C_G(g)^0$, which is also the unique $2$-Sylow$^0$ in $C_G(g)$.

Consider involutions $i\in S$ and $i'\in S'$. Then $i$ and $i'$  invert the element $m=ii'$. Let $S''$ be the unique $2$-Sylow$^0$ in $C_G(m)$. Then $i$ and $i'$ normalize $C_G(m)$, whence also $S''$. Now either $S''\not=S$ or $S''\not=S'$, and we may assume $S''\not=S$. Then $i\notin C_G(S'')$, as otherwise $S''\le C_G(i)$ would be the unique $2$-Sylow$^0$ in $C_G(i)$, contradicting $S\not=S''$.

Suppose $C_{C_G(S'')}(i)$ is infinite. Then it contains a torus $T$, and there is a full torus $A''\ge S''T$. As all full tori in $C_G(i)^0$ contain $i$, there is one, say $A^*$, which also contains $T$. But $A''\cap A^*\ge T$ so $S''\le A^*$, contradicting $i\notin C_G(S'')$. As $i$ normalizes $C_G(S'')$ and centralizes only finitely many points, it inverts $C_G(S'')^0=C_G(S'')$. It follows that $C_G(S'')$ is abelian, and $C_G(S'')=C_G(A'')=A''$; moreover all involutions in $N_G(A'')\setminus A''$ invert $A''$ and are in the same coset $iA''$ modulo $A''$; conversely all elements of the coset $iA''$ are involutions inverting $A''$. Note that $m\in C_G(S'')=A''$, and $i\in N_G(S'')=N_G(A'')$; if $j,k\in A''$ are two involutions, they commute with $i$ and invert $A$ (the full torus containing $S$). So $ik\in C_G(A)\cap A''=A\cap A''=\{1\}$, and a full torus contains a unique involution. For any non-trivial $g\in G$ the centralizer $C_G(g)$ contains a unique $2$-Sylow$^0$, and hence a unique full torus, which must contain $g$. It follows that $C_G(g)^0$ is this full torus. So the full tori of $G$, which
coincide with the connected centralizers of nontrivial elements, are disjoint and cover $G$.

Let $N$ be a minimal normal definable subgroup of $G$. Then $N$ cannot be abelian, so $N$ is simple by connectedness of $G$. 
If $N$ does not contain an involution, then any involution $i$ of $G$ gives rise to an involutive automorphism of $N$. If $F=C_N(i)$ is the subgroup of fixed points and $I$ the set of points inverted by $i$, then $N$ decomposes uniquely as $F\cdot I$ and all conjugates of $F$ intersect $I$ trivially by \cite{Po18}. Hence $F$ is connected, and a torus of $N$. But the conjugates of $F$ cover $N$, so $I=\{1\}$ and $N=F$, a contradiction. It follows that $N$ contains involutions, and we may replace $G$ by~$N$.

We obtain a geometry on $G$ whose points are the involutions, and whose lines are the cosets $iA$ for a full torus $A$ and involution $i\in N_{G}(A)$, i.e.\ sets $\ell(j)=N_{G}(A(j))[2]\setminus\{j\}$, where $A=A(j)$ is the unique full torus containing $j$ (recall that $G[2]$ is the set of involutions of $G$). Then $j$ is the unique involution commuting with all the involutions of $\ell(j)$. Any two distinct lines $\ell(i)$ and $\ell(j)$ have empty intersection, or intersect in a coset of $A(i)\cap A(j)=\{1\}$ (note that the corresponding $2$-Sylows must be distinct), i.e.\ in a single point. Any two distinct points $i,j$ lie on at least one common line $\ell(k)$, where $A(k)=C_G^0(ij)$. There is a polarity which associates an involution $i$ to a line $\ell(i)$ and vice versa, and which preserves incidence. It follows that any two lines intersect in a unique point, and any two points lie on a unique line. Moreover, there are no isotropic points: $i\notin\ell(i)$ for all $i$.
By Bachman's theorem \cite{Bach}, $G$ is a group of linear transformations of a $3$-dimensional vector space over an interpretable field $K$ which preserve a  symmetric bilinear form without isotropic vectors. But this contradicts stability.

It follows that $S$ is central. But then $G/\dc(S)$ has finite Sylow $2$-subgroups, which must be trivial by \cite{BBC}. Thus $S$ is the connected central Sylow $2$-subgroup of $G$.\end{proof}

Deloro and Wiscons have obtained Theorem \ref{t:8} as a corollary of a more general theorem on the $2$-structure of a connected group of finite Morley rank \cite[Corollary B2]{DW19}. 

\section{Characteristic $\not=2$}\label{five}
Instead of working with a split sharply $2$-transitive group, we work in the equivalent category of near-fields, see Fact \ref{f:nearfield}. Recall that addition is commutative in a near-field  by Fact \ref{fact}(\ref{commute}).
\begin{definition} The {\em kernel} $\ker(K)$ of a near-field $K$ is the set of elements with respect to which multiplication is right distributive:
\[\ker(K)=\{x\in K:\forall\,y,z\in K\,(y+z)x=yx+zx\}.\]\end{definition}
\begin{remark} The prime field of a near-field is contained in the kernel.\end{remark}
\begin{proof} Consider $n=1+\cdots+1$, and $y,z\in K$. Then by left distributivity
\[(y+z)n=(y+z)+\cdots+(y+z)=(y+\cdots+y)+(z+\cdots+z)=yn+zn.\]
In characteristic $0$ we obtain for $m=1+\cdots+1\not=0$\,:
\[\begin{aligned}(y+z)m^{-1}n&=(ym^{-1}m+zm^{-1}m)m^{-1}n=(ym^{-1}+zm^{-1})mm^{-1}n\\
&=(ym^{-1}+zm^{-1})n=ym^{-1}n+zm^{-1}n.\qedhere\end{aligned}\]
\end{proof}

Note that this does not imply that the prime field is in the centre of $K$, even if the kernel is finite (and $K$ connected), as conjugation is not an automorphism of $K$ and need not stabilize the kernel.

\begin{theorem}\label{nearfield} An infinite near-field $K$ of finite Morley rank in characteristic $\not=2$ is an algebraically closed field.\end{theorem}
\begin{proof} If the kernel is infinite (in particular if char$(K)=0$ or $Z(K^\times)$ is infinite), this follows from \cite{BN} or \cite{CGNV}.

In characteristic $p>0$, note first that $K$ is additively connected, as for any additive proper subgroup $H$ of finite index the intersection $\bigcap_{x\in K^\times}xH$ is trivial, but equals a finite subintersection, and hence is of finite index, a contradiction. Therefore there is a unique type of maximal Morley rank, so $K^\times$ is multiplicatively connected as well.

Let $A$ be a definable connected infinite abelian multiplicative subgroup, and $M_0$ an $A$-minimal additive subgroup. By Zilber's Field Theorem $M_0$ is additively isomorphic to the additive group of an algebraically closed field $K_0$, and $A$ embeds multiplicatively into $K_0^\times$. In fact, for any $e_0\in M_0\setminus\{0\}$ and $a_1,\ldots,a_n\in A$ such that $a_1e_0+\cdots+a_ne_0=0$, we have by left distributivity that $a_1^{e_0}+\cdots+a_n^{e_0}=0$, so the addition induced on $A$ by $K_0$ is the one inherited from $K$-addition on $A^{e_0}$. So we might replace $A$ by $A^{e_0}$, $M_0$ by $e_0^{-1}M_0$  and $e_0$ by $1$. Then $A\subseteq M_0=K_0^+$, and field multiplication on $K_0$ is induced from $K$ on $A\times K_0$, but does not necessarily agree with multiplication from $K$ if the left factor is in $K_0\setminus A$. In particular, $A$ is a good torus by \cite{Wa}. By Theorem \ref{t:8} the centre $Z(K^\times)$ contains the Sylow $2$-subgroup, which is infinite in characteristic $\not=2$. We finish by the first paragraph.
\end{proof}

\begin{corollary} A split sharply $2$-transitive permutation group of finite Morley rank and characteristic $\not=2$ is standard. In particular a sharply $2$-transitive permutation group of finite Morley rank and characteristic $3$ is standard.\end{corollary}
\begin{proof} By \cite{KW} a sharply $2$-transitive permutation group of characteristic $3$ splits. Now use Fact \ref{f:nearfield} and Theorem \ref{nearfield} to finish.\end{proof}


\begin{thebibliography}{99}

\bibitem{ABC}
T. Alt\i nel, A. Borovik, and G. Cherlin,
\emph{Simple Groups of Finite Morley Rank}.
Math.\ Surv.\ Monogr.\ 145, AMS, 2008.

\bibitem{AB08}
T. Alt\i nel and J. Burdges,
On analogies between algebraic groups and groups of finite {M}orley rank.
\emph{J. Lond.\ Math.\ Soc.} \textbf{78} (2008), no.~1, 213--232.

\bibitem{AAT} M. Amelio, S. Andr\'{e}, K. Tent, Non-split sharply 2-transitive groups of odd positive characteristic (2023). Int.\ Math.\ Res.\ Not.\ IMRN (2025), no.\ 19, article no. rnaf294.
		
\bibitem{AG} S.	Andr\'{e}, V. Guirardel, Finitely generated simple sharply 2-transitive groups. \emph{Comp. Math.} \textbf{160} (2024), no.~8, 1941--1957. 
		
\bibitem{AT} S.	Andr\'{e}, K. Tent, 
Simple sharply 2-transitive groups.
\emph{Trans.\ Amer.\ Math.\ Soc.} \textbf{376} (2023), no.~6, 3965--3993.

\bibitem{Bach}
F. Bachmann,
\emph{Aufbau der Geometrie aus dem Spiegelungsbegriff}.
Springer Verlag, 1973.

\bibitem{BBC}
A. Borovik, J. Burdges, and G. Cherlin,
Involutions in groups of finite {M}orley rank of degenerate type.
\emph{Sel.\ Math.\ New Ser.} \textbf{13} (2007), no.~1, 1--22.

\bibitem{BN} 
A. Borovik and A. Nesin,
\emph{Groups of finite Morley rank}.
Oxford Logic Guides 26, OUP, 1994.

\bibitem{Ch05}
G. Cherlin,
Good tori in groups of finite {M}orley rank.
\emph{J. Group Theory} \textbf{8} (2005), no.~5, 613--622.

\bibitem{CGNV}
G. Cherlin, T. Grundh\"ofer, A. Nesin and H. V\"olklein,
Sharply Transitive Linear Groups over Algebraically Closed Fields.
\emph{Proc.\ Amer.\ Math.\ Soc.} \textbf{111} (1991), no.~2, 541--550.

\bibitem{CT2} T. Clausen and K. Tent, Mock hyperbolic reflection spaces and Frobenius groups of finite Morley rank. \emph{Model Theory} \textbf{2} (2023), no.~2, 137--175.

\bibitem{DW19}
A. Deloro and J. Wiscons, 
The geometry of involutions in ranked groups with a TI-subgroup.
 \emph{Bull.\ Lond.\ Math.\ Soc.} \textbf{52} (2020), no.~3, 411--428.

\bibitem{Fr01}
F. G. Frobenius,
{\"U}ber aufl{\"o}sbare {G}ruppen IV.
\emph{Berl.\ Ber.} (1901), 1216--1230.

\bibitem{GG} Y. Glasner, D. Gulko, 
Non-split linear sharply 2-transitive groups. \emph{Proc.\ Am.\ Math.\ Soc.} \textbf{149} (2021), No.~6, 2305--2317.

\bibitem{GG2} Y. Glasner, D. Gulko, Sharply 2-transitive linear groups. \emph{Int.\ Math.\ Res.\ Not.} \textbf{2014} (2014), No.~10, 2691--2701.

\bibitem{Ka}
H. Karzel,
Zusammenh\"ange zwischen {F}astbereichen, scharf $2$-fach transitiven {P}ermutationsgruppen und $2$-{S}trukturen mit {R}echtecksaxiom.
\emph{Abh.\ Math. Sem.\ Univ.\ Hamburg} \textbf{32} (1968), 191--206.

\bibitem{Ker74} Kerby, William. {\em On infinite sharply multiply transitive groups}, Hamburger Mathematische Einzelschriften, Neue Folge, Heft 6. Vandenhoeck \& Ruprecht, G\"ottingen 1974.

\bibitem{KW}
W. Kerby and H. Wefelscheid,
Bemerkungen \"uber {F}astbereiche und scharf $2$-fach transi\-tive {G}ruppen.
\emph{Abh.\ Math. Sem.\ Univ.\ Hamburg} \textbf{32} (1968), 191--206.

\bibitem{Kh} E. Khukhro and V. Mazurov. {\em The Kourovka Notebook}. \url{https://kourovka-notebook.org/}

\bibitem{Neu}
B. H. Neumann,
On the commutativity of addition.
\emph{J. Lond.\ Math.\ Soc.} \textbf{15} (1940), 203--208.

\bibitem{Po18}
B. Poizat,
Milieu et sym\'etrie, une \'etude de la convexit\'e dans les groupes sans involutions.
\emph{J. Alg.} \textbf{497} (2018), 143--163.

\bibitem{RST}
E. Rips, Y. Segev, and K. Tent,
A sharply 2-transitive group without a non-trivial abelian normal subgroup. \emph{J.\ Eur.\ Math.\ Soc.} \textbf{19} (2017), no.~10, 2895--2910.
	
\bibitem{RT} E. Rips, K. Tent,
Sharply 2-transitive groups of characteristic 0. \emph{J.\ Reine Angew.\ Math.} \textbf{750} (2019), 227--238.

\bibitem{ST} M. Scherff, K. Tent,
Addendum to: Sharply 2-transitive groups of characteristic 0.  \emph{J.\ Reine Angew.\ Math.} \textbf{750}  (2019), 239--240.

\bibitem{tao}
T. Tao,
The theorems of {F}robenius and {S}uzuki on finite groups.
Blog notes 04/12/2013,
\url{terrytao.wordpress.com/2013/04/12/}

\bibitem{TZ}K. Tent, and M. Ziegler, 
Sharply 2-transitive groups. \emph{Adv.\ Geom.} \textbf{16} (2016), no.~1, 131--134.

\bibitem{Th60}
J. G. Thompson,
Normal $p$-complements for finite groups.
\emph{Math.\ Zeitschr.} \textbf{72} (1960), 332--354.

\bibitem{tuek}
Seyfi T\"urkelli, Splitting of Sharply $2$-Transitive Groups of Characteristic $3$. {\em Turk.\ J. Math.} \textbf{28}(3) (2004), 295--298.

\bibitem{Wa}
F. O. Wagner, 
Fields of finite {M}orley rank.
\emph{J. Symb.\ Logic} \textbf{66} (2001), no.~2, 703--706.

\bibitem{W} F. Wagner,  Some remarks on sharply 2-transitive groups and near-domains. (2022)  arXiv:2202.13740v3.

\end{thebibliography}
\end{document}